\documentclass[12pt]{amsart}

\usepackage{aliascnt}
\usepackage{hyperref}
\usepackage{amsfonts}
\usepackage{mathrsfs}
\usepackage{amsmath}
\usepackage{amsmath,graphics}
\usepackage{amssymb}
\usepackage{indentfirst,latexsym,bm,amsmath,pstricks,amssymb,
amsthm,graphicx,fancyhdr,float,color}
\usepackage[all]{xy}
\usepackage{amsthm}
\usepackage{amscd}
\usepackage[all]{hypcap}

\newtheorem{theorem}{Theorem}[section]

\newaliascnt{lemma}{theorem}
\newtheorem{lemma}[lemma]{Lemma}
\aliascntresetthe{lemma}

\newaliascnt{conjecture}{theorem}

\aliascntresetthe{conjecture}

\newaliascnt{proposition}{theorem}
\newtheorem{proposition}[proposition]{Proposition}
\aliascntresetthe{proposition}

\newaliascnt{corollary}{theorem}
\newtheorem{corollary}[corollary]{Corollary}
\aliascntresetthe{corollary}

\newaliascnt{problem}{theorem}

\aliascntresetthe{problem}

\newaliascnt{claim}{theorem}

\aliascntresetthe{claim}

\theoremstyle{definition}

\newaliascnt{definition}{theorem}

\aliascntresetthe{definition}

\newaliascnt{example}{theorem}

\aliascntresetthe{example}

\theoremstyle{remark}

\newaliascnt{remark}{theorem}
\newtheorem{remark}[remark]{Remark}
\aliascntresetthe{remark}

\newaliascnt{remarks}{theorem}

\aliascntresetthe{remarks}

\numberwithin{equation}{section}

\numberwithin{figure}{section}

\def\wt{\widetilde}

\def\lra{\longrightarrow}

\def\({$($}
\def\){$)$}

\def\bbp{\mathbb P}

\def\call{\mathcal L}

\def\rank{\text{{\rm rank\,}}}

\begin{document}

\title{Singular Fibers and Kodaira Dimensions}

\author{Xin Lu}
\address{Department of Mathematics,
    East China Normal University, Dongchuan RD 500,
    Shanghai 200241, P. R. of China}
\curraddr{Institut f\"ur Mathematik, Universit\"at Mainz, Mainz,
55099, Germany} \email{x.lu@uni-mainz.de}

\author{Sheng-Li Tan}
\address{Department of Mathematics, East China Normal University,
    Dongchuan RD 500, Shanghai 200241, P. R. of China}
\email{sltan@math.ecnu.edu.cn}

\author{Kang Zuo}
\address{Institut f\"ur Mathematik, Universit\"at Mainz,
Mainz, 55099, Germany}
\email{zuok@uni-mainz.de}

\thanks{This work is supported by SFB/Transregio 45 Periods,
 Moduli Spaces and Arithmetic of Algebraic Varieties of DFG,
 by NSF of China and by the Science Foundation of Shanghai
 (No. 13DZ2260400).}






\phantomsection
\begin{abstract}
\addcontentsline{toc}{section}{Abstract}  Let $f:\,X \to \bbp^1$ be
a non-isotrivial semi-stable family of varieties of dimension $m$
over $\bbp^1$ with $s$ singular fibers. Assume that the smooth
fibers $F$ are minimal, i.e., their canonical line bundles are
semiample. Then $\kappa(X)\leq \kappa(F)+1$. If
$\kappa(X)=\kappa(F)+1$, then $s>\frac{4}m+2$. If $\kappa(X)\geq 0$,
then $s\geq\frac{4}m+2$. In particular, if $m=1$, $s=6$ and
$\kappa(X)=0$, then the family $f$ is Teichm\"uller.
\end{abstract}

\maketitle

\section{Introduction}
We always work over the complex number $\mathbb C$. Let $f:\,S \to
\mathbb P^1$ be a nontrivial fibration of semi-stable curves of
genus $g\geq 1$. It is a classical problem to determine the lower
bound for the number $s$ of singular fibers in the fibration $f$,
see  \cite{bea-81,tan-95,ttz-05,tu-07,zamora-12,glt-13,ltxz-16}. In
\cite{bea-81}, Beauville first proved that $s\geq 4$ and conjectured
that $s\geq 5$ when $g\geq2$. In \cite{tan-95}, the second author
confirmed Beauville's conjecture. Later, Tu, Zamora and the second
author proved in \cite{ttz-05} that $s\geq 6$ if $S$ has
non-negative Kodaira dimension. It is conjectured that $s\geq 7$ if
$S$ is of general type. The first purpose of this note is to confirm
this conjecture.

\begin{theorem}\label{thm-main}
Let $f:\,S \to \mathbb P^1$ be a nontrivial semi-stable
fibration $f:\,S\to\bbp^1$
of curves of genus $g\geq 2$ over $\bbp^1$ with $s$ singular fibers.
If $S$ is of general type, then $s\geq 7.$
\end{theorem}

 This conjecture has been verified for $g\leq 5$
 (\cite{ttz-05,zamora-12} or $g\geq 58$ (\cite{TTY},
 unpublished) by using the {\it strict canonical class
 inequality} established by the second author \cite{tan-95}.
 Recently, the authors in \cite{msz-16}
have also proved this conjecture under the condition that the family
is birationally equivalent to a pencil of curves with only simple
base points on the minimal model of $S$.

We can find in  \cite{ttz-05} the examples of surfaces
of general type admitting a semi-stable fibration over
$\bbp^1$ with $7$ singular fibers.

It is an interesting phenomenon that when the number of singular
fibers is minimal, the family is of very interesting arithmetic and
geometric properties. When $g=1$ and $s=4$, Beauville \cite{bea-82}
proved that the family of curves must be {\it modular}, and there
are exactly 6 such families. In \cite{STZ04}, the authors prove that
for a non-isotrivial family of semi-stable K3 surfaces
$f:X\to\mathbb P^1$ on a Calabi-Yau manifold $X$, then $s\geq 4$ and
if $s=4$, the family is {\it modular}.

\begin{theorem}\label{20161018a} As in \autoref{thm-main},
if the Kodaira dimension
of $S$ is zero and $s=6$, then the family must be Teichm\"uller
and $\omega^2_{S/\mathbb P^1}=6g-6$.
\end{theorem}

Here a family of curves is said to be Teichm\"uller, if
up to a suitable finite \'etale cover of $S$,
it comes from a Techm\"uller curve.

For each type of surfaces of Kodaira dimension zero,
Tu \cite{tu-07}
 has constructed an example with a semi-stable family of
 curves over $\bbp^1$ admitting exactly $6$ singular fibers.

When the Kodaira dimension of the surface is $1$, the
minimal number $s$ should be 6 or 7.
We have not found examples with $s=6$, and we tend to
believe that there are no such examples.
On the other hand, we give more precise description
of such surfaces.

\begin{theorem}\label{thm-3-1}
    With the notation as in \autoref{thm-main}.
     Suppose the Kodaira dimension of $S$ is $1$ and $s=6$. Then
    $S$ is simply connected,
    $p_g(S)=q(S)=0$,  the canonical elliptic fibration
    on $S$ admits exactly two multiple fibers,
    one of the multiplicities is $2$,  and the second one is $n=3$ or $5$.
        \begin{enumerate}
            \item[(1)] If $n=3$, then $6g-5\leq
            \omega_{S/\bbp^1}^2\leq 6g-3$.
            \item[(2)] If $n=5$, then $\omega_{S/\bbp^1}^2= 6g-3$.
        \end{enumerate}
\end{theorem}

When the stability assumption is dropped, then it is only
known that $s\geq 3$ for any non-isotrivial fibration of curves
    over $\bbp^1$, even if we require that two of the
    singular fibers be semi-stable \cite{bea-81,glt-13}.

Our method works also for the higher dimensional cases.
\begin{theorem}\label{thm-high-dim}
    Let $f:\,X \to \bbp^1$ be a non-isotrivial semi-stable
    family of varieties of dimension $m$ over $\bbp^1$ with
    $s$ singular fibers.
Assume that the smooth fibers $F$ are minimal, i.e., their
canonical line bundles are semiample.
Then $\kappa(X)\leq \kappa(F)+1$.
\begin{enumerate}
\item[{\rm 1)}] If $\kappa(X)\geq 0$, then $s\geq\frac{4}m+2$.
In particular, $s\geq 6$ when $m=1$, and $s\geq 4$
    when $m=2$ or $3$.
\item[{\rm 2)}] If $\kappa(X)=\kappa(F)+1$, then $s>\frac{4}m+2$.
In particular, $s\geq 7$ when $m=1$, $s\geq 5$ when $m=2$,
 and $s\geq 4$ when $m=3$ or $4$.
\end{enumerate}
\end{theorem}

Note that when $X$ is of general type, then $F$ must be also of
general type and the equality $\kappa(X)=\kappa(F)+1$ holds.
Hence the lower bound $s>\frac{4}m+2$ holds in this case.

\vspace{2mm}
This note is organized as follows.  Theorems \ref{thm-main},
\ref{20161018a} and \ref{thm-3-1} are proved in \autoref{sec-2},
and Theorem \ref{thm-high-dim} is proved in \autoref{sec-3}.

\section{Variations of the Hodge
structures}\label{sec-2} In this section, we would like to prove
Theorems \ref{thm-main}-\ref{thm-3-1}. The main technique is based
on the variation of the Hodge structures attached to a semi-stable
family of curves, especially to a techm\"uller family.

\subsection{Preliminaries}
In this subsection, we give a brief recall about the
Techm\"uller curve and the associated variation of the
Hodge structures,
and derive some inequalities.
For more details, we refer to \cite{mol-06,mol-13,EKZ-14}.

Let ${\mathcal M_g}$ be the moduli space of smooth projective curves
of genus $g$, and $\Omega {\mathcal M_g}\to {\mathcal M_g}$ the
bundle of pairs $(F,\omega)$, where $\omega\neq 0$ is a holomorphic
one-form on $F \in {\mathcal M_g}$. Let $\Omega {\mathcal
M_g}(m_1,\cdots,m_k) \subseteq \Omega {\mathcal M_g}\to {\mathcal
M_g}$ be the stratum of pairs $(F,\omega)$ such that $\omega$ admits
exactly $k$ distinct zeros of order $m_1,\cdots, m_k$ respectively.
There is a natural action of $SL_2(\mathbb R)$ on each stratum
$\Omega {\mathcal M_g}(m_1,\cdots,m_k)$. Each orbit projects to a
complex geodesics in ${\mathcal M_g}$. When the projection of such
an orbit is closed, it gives a so-called {\it Teichm\"uller curve}.
After a suitable unramified cover and compactification of a given
Teichm\"uller curve, one gets a universe family $f:\,S \to B$, which
is a semi-stable family of curves of genus $g$. Moreover, there
exist disjoint sections $D_1,\cdots D_k$ of $f$ such that the
restriction $\big(\sum\limits_{i=1}^k m_iD_i\big)\big|_F$ to each
fiber $F$ is just the zero locus of $\omega$.

Denote by $s$ the number of singular fibers contained in $f$.
Then the Hodge bundle $f_*\omega_{S/B}$ for a Teichm\"uller
curve contains
a line subbundle $\mathcal L\subseteq f_*\omega_{S/B}$ with
maximal slope:
\begin{equation}\label{eqn-3-1}
2\deg(\mathcal L)=2g(B)-2+s.
\end{equation}

Consider the logarithmic Higgs bundle $(f_*\omega_{S/B}
\oplus R^1f_*\mathcal O_{S},\,\theta)$ associated to the
fibration $f$,
which corresponds to the weight-one local system
$R^1f_*\mathbb{Q}_{S^0}$;
here $f:\,S^0 \to B^0$ is the smooth part of $f$.
The Higgs field $\theta$ is simply the edge morphism
$$f_*\omega_{S/B}\cong f_*\Omega^1_{S/B}(\log\Upsilon)
\lra R^1f_*\mathcal O_{S} \otimes \Omega^1_{B}(\log\Delta)$$
of the tautological sequence
$$0\lra f^*\Omega^1_{B}(\log\Delta)\lra \Omega^1_{S}(\log\Upsilon)
\lra \Omega^1_{S/B}(\log\Upsilon) \lra 0,$$ where $\Upsilon \to
\Delta$ is denoted to be the singular locus of $f$. By
\cite{viehweg-zuo-03}, the existence of a line subbundle $\mathcal L
\subseteq f_*\omega_{S/B}$ with maximal slope is equivalent to the
existence of a rank two Higgs subbundle $(\mathcal L\oplus \mathcal
L^{-1},\,\theta)$ with maximal Higgs field contained in the
logarithmic Higgs bundle $(f_*\omega_{S/B} \oplus R^1f_*\mathcal
O_{S},\,\theta)$ associated to the fibration $f$.

Conversely, one has the following theorem, which is due to
M\"oller \cite{mol-06}.
\begin{theorem}\label{thm-mol}
    Let $f:\,S \to B$ be a semi-stable fibration of curves of
    genus $g\geq 2$ over a smooth projective curve with $s$
    singular fibers.
    Suppose that there exists a line subbundle
    $\mathcal L\subseteq f_*\omega_{S/B}$ satisfying the
    equality \eqref{eqn-3-1} above.
    Then the family $f$ comes from a Techm\"uller curve;
    that is, the induced map $B^0\to {\mathcal M_g}$ is a finite
    unramified cover of a Techm\"uller curve.
    Here $f:\,S^0 \to B^0$ is the smooth part of $f$.
\end{theorem}

Since the relative canonical sheaf of a fibration of curves
over a Techm\"uller curve has a very special form
(see \eqref{eqn-3-2} below),
we can derive the following upper bound on $\omega_{S/B}^2$.

\begin{proposition}\label{prop-3-1}
    Let $f:\,S \to B$ be a semi-stable fibration of curves
    as in \autoref{thm-mol},
    and assume also that there exists a line subbundle
    $\mathcal L\subseteq f_*\omega_{S/B}$ with the
    equality \eqref{eqn-3-1}.
    Then
    \begin{equation}\label{eqn-3-3}
    \omega_{S/B}^2\leq \frac{3}{2}(g-1)\big(2g(B)-2+s\big).
    \end{equation}
\end{proposition}
\begin{proof}
    By \autoref{thm-mol}, the induced map $B^0\to {\mathcal M_g}$ is
    finite unramified covering of a Techm\"uller curve.
    Hence after a suitable unramified base change,
    there exist disjoint sections $D_1,\cdots D_k$ of $f$
    such that
    the relative canonical sheaf $\omega_{S/B}$ has the
    form (cf. \cite{EKZ-14})
    \begin{equation}\label{eqn-3-2}
    \omega_{S/B}\cong f^*\mathcal L \otimes
    \mathcal O_S\Big(\sum_{i=1}^k m_iD_i\Big),
    \end{equation}
    where $\mathcal L \subseteq f_*\omega_{S/B}$ is the line
    subbundle satisfying the equality \eqref{eqn-3-1}.
    Note that the inequality \eqref{eqn-3-3} is invariant
    under any finite unramified base change.
    Thus we may assume that $\omega_{S/B}$ already has the
    form as above.

    As $D_i$'s are disjoint sections, it follows that
    $D_i\cdot D_j=0$ for $i\neq j$, and that
    $$(\omega_{S/B}+D_i)\cdot D_i=0,\quad \forall~1\leq i\leq k.$$
    Combining these with \eqref{eqn-3-2}, one gets that
    $D_i^2=-\frac{1}{m_i+1}\cdot \deg\call$.
    Note also that $\sum\limits_{i=1}^km_i=\deg\omega_F=2g-2$,
    where $\omega_F$ is the canonical sheaf on a general
    fiber of $f$.
    Hence by \eqref{eqn-3-2} again, we obtain that
    \begin{align*}
    \omega_{S/B}^2&=4(g-1)\cdot \deg\call+
    \sum_{i=1}^k m_i^2D_i^2\\
    &=\Big(4(g-1)-
    \sum_{i=1}^k \frac{m_i^2}{m_i+1}\Big)\deg\call.
    \end{align*}
    As $\sum\limits_{i=1}^km_i=2g-2$, one gets easily that
    $$\sum_{i=1}^k \frac{m_i^2}{m_i+1}
    \geq \sum_{i=1}^k \frac{m_i}{2}= g-1.$$
    Therefore,
    $$\omega_{S/B}^2 \leq 3(g-1)\cdot
    \deg \call=\frac{3(g-1)\big(2g(B)-2+s\big)}{2}.$$
    This completes the proof.
\end{proof}

In the case when $f:\,S \to \bbp^1$ is a semi-stable
fibration of curves of genus $g\geq 2$ over $\bbp^1$
with $s=6$ singular fibers,
we have the following easy criterion when $f$ comes
from a Techm\"uller curve.
\begin{lemma}\label{cor-3-1}
    Let $f:\,S \to \bbp^1$ be a semi-stable fibration
    of curves of genus $g\geq 2$ over $\bbp^1$ with $s=6$
    singular fibers.
    If the geometric genus $p_g(S):=\dim H^0(S,\omega_S)>0$,
    then there exists a line subbundle
    $\mathcal L\subseteq f_*\omega_{S/B}$ satisfying the
    equality \eqref{eqn-3-1}, and hence $f$ comes from a
    Techm\"uller curve.
\end{lemma}
\begin{proof}
    As a locally free sheaf on $\bbp^1$,  the direct image
    sheaf $f_*\omega_{S/\bbp^1}$ is isomorphic to
    a direct sum of invertible sheaves:
    $$f_*\omega_{S/\bbp^1} \cong \bigoplus_{i=1}^{g}
    \mathcal O_{\bbp^1}(d_i).$$
    Note that $d_i\geq 0$ due to the semi-positivity of the
    direct image sheaf $f_*\omega_{S/\bbp^1}$ (cf. \cite{fujita-78}),
    and that $d_i\leq \frac{1}{2}\big(2g(B)-2+s\big)=2$
    due to the Arakelov type inequality (cf. \cite{viehweg-zuo-03}).
    Without loss of generality, we assume
    that $0\leq d_1\leq \cdots \leq d_g\leq 2$.
    By \cite[Theorem\,3.1]{fujita-78}, we obtain that
    $$d_1=\cdots =d_{q(S)}=0,\qquad \text{and}\quad d_i>0
    \quad \forall~ i\geq q(S)+1,$$
    where $q(S):=\dim H^1(S,\omega_S)$ is the irregularity of $S$.
    Hence
    $$\deg f_*\omega_{S/\bbp^1} =\sum_{i=q(S)+1}^{g} d_i.$$
    On the other hand, it is well-known that
    $$\deg f_*\omega_{S/\bbp^1}=\chi(\mathcal \omega_S)-(g-1)
    \big(g(\bbp^1)-1\big)=g+p_g(S)-q(S).$$
    Therefore, $d_g=2$ once $p_g(S)>0$. In other word, the
    line subbundle $\mathcal O_{\bbp^1}(d_g)
    \subseteq f_*\omega_{S/\bbp^1}$
    satisfies the equality \eqref{eqn-3-1}.
\end{proof}

\begin{corollary}\label{cor-3-2}
    Let $f:\,S \to \bbp^1$ be a semi-stable fibration of curves
    of genus $g\geq 2$ over $\bbp^1$ with $s=6$ singular fibers.
    If the geometric genus $p_g(S)>0$, then
    $f$ comes from a Techm\"uller curve and
    \begin{equation}\label{eqn-3-4}
    \omega_{S/\bbp^1}^2\leq 6(g-1).
    \end{equation}
\end{corollary}
\begin{proof}
    This is a combination of \autoref{cor-3-1} and \autoref{prop-3-1}.
\end{proof}

\subsection{Proof of \autoref{thm-main}}\label{sec-2-2}
By \cite[Theorem\,0.1]{ttz-05}, $s\geq 6$ if $S$ is of general type
(actually, the inequality $s\geq 6$ holds once the Kodaira
dimension of $S$ is non-negative).
To complete the proof, it suffices to deduce a contradiction
if $s=6$.

Since $S$ is of general type,
 we may assume that $g\geq 5$ by \cite[Theorem\,0.1(2)]{ttz-05},
and according to \cite[Theorem\,0.2]{ttz-05} one has
\begin{equation}\label{eqn-2-21}
    \omega_{S/\bbp^1}^2\geq 6g-6+\frac12\left(\omega_X^2+
    \sqrt{\omega_X^2}\,\sqrt{\omega_X^2+8g-8}\right),
\end{equation}
where $X$ is the minimal model of $S$.
Hence we may assume that $p_g(S)=0$ by  \autoref{cor-3-2}.
It then follows that
$$\deg f_*\omega_{S/\bbp^1}=\chi(\mathcal \omega_S)-(g-1)
\big(g(\bbp^1)-1\big)=g-q(S).$$
Let $\delta_f$ be the number of nodes contained in the
fibers of $f$. Then by Noether's formula, one has
$$\delta_f=12 \deg f_*\omega_{S/\bbp^1}-\omega_{S/\bbp^1}^2
=12\big(g-q(S)\big)-\omega_{S/\bbp^1}^2.$$
According to \cite{tan-95}, for any integer $e\geq 2$, we
have the following inequality
$$\begin{aligned}
\omega_{S/\bbp^1}^2 &~\leq (2g-2)\Big(2g\big(\bbp^1\big)-2+
\frac{(e-1)s}{e}\Big)+\frac{3\delta_f}{e^2},\\
&~=(2g-2)\Big(4-\frac{6}{e}\Big)+\frac{3\left(12\big(g-q(S)\big)
-\omega_{S/\bbp^1}^2\right)}{e^2}.
\end{aligned}$$
Hence
$$\omega_{S/\bbp^1}^2 \leq \frac{e^2}{e^2+3}(2g-2)
\Big(4-\frac{6}{e}\Big)+\frac{36\big(g-q(S)\big)}{e^2+3}.$$
Taking $e=3$, one obtains
\begin{equation}\label{eqn-3-5}
\omega_{S/\bbp^1}^2 \leq 6g-3-3q(S) \leq 6g-3.
\end{equation}
Combining this with \eqref{eqn-2-21}, one obtains that
$$\begin{aligned}
&~\sqrt{\omega_X^2}\,\sqrt{\omega_X^2+8g-8} \leq 6-\omega_X^2,\\
\Longrightarrow ~&~\omega_X^2 \,\big(\omega_X^2+8g-8\big)
\leq (6-\omega_X^2)^2,\\
\Longrightarrow ~&~ \omega_X^2\leq \frac{9}{2g+1} <1,
\qquad \text{since~}g\geq 5.
\end{aligned}$$
This gives a contradiction.
\qed

\subsection{Proof of \autoref{20161018a}}\label{sec-2-3}
    If $S$ is either an abelian surface or a {\rm K3} surface,
    then $p_g(S)>0$,
    and hence the conclusion follows directly
    from \autoref{cor-3-2} and \cite[Theorem\,0.2]{ttz-05}.

    In the remaining cases, $S$ must be either an Eriques
    surface or a bielliptic surface
    according to the classification of surfaces with Kodaira
    dimension equal to zero.
    Let $K_X$ be the canonical divisor on the minimal model
    $X$ of $S$.
    Then there exists an $n>1$ such that $nK_X\equiv 0$.
    Hence one can construct a finite \'etale cover $\pi:\,\wt S \to S$
    such that $p_g(\wt S)>0$ and that the Kodaira dimension of $\wt S$
    is still zero.
    Moreover $\tilde f:=f\circ \pi:\,\wt S \to \bbp^1$ is still a
    semi-stable fibration with 6 singular fibers by
    \cite[Lemma\,3]{bea-81}.
    $$\xymatrix{
        \wt S \ar[rr]^-{\pi} \ar[dr]_-{\tilde f} && S \ar[dl]^-{f}\\
        &\bbp^1&}$$
    Since $\pi$ is finite \'etale,
    \begin{equation}\label{eqn-3-7}
    \omega_{\wt S/\bbp^1}^2=\deg(\pi)\cdot \omega_{S/\bbp^1}^2,\qquad
    \tilde g-1=\deg(\pi)\cdot (g-1),
    \end{equation}
    where $\tilde g$ is the genus of a general fiber of $\tilde f$.
    By construction, $\wt S$ is either an abelian surface or
    a {\rm K3} surface, so $p_g(\tilde S)=1$. Hence the family $\tilde f$ comes from
    a Techm\"uller curve
    and $\omega_{\wt S/\bbp^1}^2=6(\tilde g-1)$ by the above argument.
    Therefore, the family $f$ is Techm\"uller. Moreover,
    Together with \eqref{eqn-3-7}, we obtain $\omega_{S/\bbp^1}^2
    = 6g-6$ as required.
\qed

\subsection{Proof of \autoref{thm-3-1}}
    Because the Kodaira dimension of $S$ is $1$,
    by \cite[Theorem\,0.2]{ttz-05} we have
    \begin{equation}\label{eqn-2-22}
        \omega_{S/\bbp^1}^2 \geq 6g-5.
    \end{equation}
    Hence $p_g(S)=0$ by \autoref{cor-3-2}.
    Similar to the proof of \autoref{thm-main},
    the inequality \eqref{eqn-3-5} holds.
    Thus $q(S)=0$ by \eqref{eqn-2-22} and \eqref{eqn-3-5}.

    As the kodaira dimension of $S$ is $1$,
    the minimal model $X$ of $S$ admits an elliptic fibration
    $$h:\,X \lra C.$$
    Since $q(X)=q(S)=0$, it follows that $C\cong \bbp^1$.
    Let $\{n_1\Gamma_1,\cdots,n_r\Gamma_r\}$ be the set of
    multiple fibers of $h$ with $2\leq n_1\leq \cdots\leq n_r$.
    Then the canonical sheaf of $X$ is given by
    (cf. \cite[\S\,IV-5]{gh-94})
    \begin{equation}\label{eqn-3-6}
    \omega_{X} = h^*\Big(\mathcal O_{\bbp^1}(-1)\Big)\otimes
    \mathcal O_X\Big(\sum_{i=1}^{r}(n_i-1)\Gamma_i\Big).
    \end{equation}

    We claim first that $r=2$.  Indeed, it is clear that
    $r\geq 2$ by \eqref{eqn-3-6} since $\kappa(X)=1$,
    and that $r<3$, since otherwise by an unramified cover
    one can construct a new surface $\wt S$ with $p_g(\wt S)>0$.
    Moreover, similar to the proof of \autoref{20161018a},
    one shows that  $\wt S$ is still a
    semi-stably fibred over $\bbp^1$ with $6$ singular fibers.
    This is a contradiction by the above argument.

    We claim also that $n_1{\not|~}n_2$. Suppose
    $n_1$ divides $n_2$,  one can construct
    an unramified cover $S''$ over $S$,
    which is still semi-stably fibred over $\bbp^1$ with $6$
    singular fibers.
    Moreover, the minimal model of $S''$ admits an elliptic
    fibration with only one multiple fiber.
    This is again a contradiction by the above argument.

    Let $F$ be a general fiber of $f$ and $F_0$ its image in $X$.
    Let $\Gamma$ be a general fiber of $h$ and $d=\gcd(n_1,n_2)$.
    Then there exist $m_1, m_2\in \mathbb Z$ such that
    $m_1n_1+m_2n_2=d$.
    Let $\Gamma_0=m_2\Gamma_1+m_1\Gamma_2$.
    Then numerically,
    $$
    \Gamma_0 = \frac{m_2}{n_1} \cdot n_1\Gamma_1+
    \frac{m_1}{n_2}\cdot n_2\Gamma_2
    \sim_{\rm num} \left(\frac{m_2}{n_1} +\frac{m_1}{n_2}\right)
    \Gamma =\frac{d}{n_1n_2}\cdot \Gamma.$$
    Moreover,  by \eqref{eqn-3-6}, one has the following numerical
    equivalence:
    $$\omega_X\sim_{\rm num} \left(1-\frac{1}{n_1}
    -\frac{1}{n_2}\right)\Gamma
    \sim_{\rm num} \left(\frac{n_1n_2}{d}-\frac{n_1}{d}
    -\frac{n_2}{d}\right) \Gamma_0.$$
    According to the proof of \cite[Theorem\,2.1]{ttz-05},
    one has
    $$\begin{aligned}
    \omega_{S/\bbp^1}^2~&\geq 6g-6+\omega_X\cdot F_0\\
    ~&=6g-6+ \left(\frac{n_1n_2}{d}-\frac{n_1}{d}
    -\frac{n_2}{d}\right) \Gamma_0\cdot F_0\\
    ~&\geq 6g-6+ \left(\frac{n_1n_2}{d}-\frac{n_1}{d}
    -\frac{n_2}{d}\right).
    \end{aligned}$$
    From $n_1{\not|~}n_2$ and \eqref{eqn-3-5},
    we see that there are only two possibilities as
    stated in \autoref{thm-3-1}.

    It remains to show that $S$ is simply connected.
    Since $\chi(\mathcal{O}_X)=1>0$,
    it follows from Noether's formula that the elliptic fibration
    $h$ admits at least one singular fiber.
    Moreover, we have shown that $h$ has exactly two
    multiply fibers whose multiplicities are coprime.
    From \cite[\S\,II.2-Theorem\,10]{moi-77},
    it follows that $X$, and hence also $S$, are both
    simply connected.
    This completes the proof.
\qed

\section{Arakelov type inequality}\label{sec-3}
In this section, we generalize our results to the
high dimension cases,
i.e., we prove \autoref{thm-high-dim}.
The technique uses the Arakelov type inequality,
which is deduced from the variation of the Hodge
structures attached to such families.

The Arakelov type inequality for the direct image of the relative
pluri-canonical sheaves goes back to Viehweg and the last author
\cite{viehweg-zuo-06,zuo-08}. This kind of inequality is generalized
in the recent work \cite{ltz-16}. The following form can be found in
\cite[Theorem\,4.4]{zuo-08} and \cite[Prop\,3.1 \&
Remark\,3.2]{ltz-16}, which is the key to our proof.
\begin{theorem}\label{thm-ltz}
    Let $f:\,X\to B$ be a semi-stable family of varieties of
    relative dimension $m\geq 1$ over a smooth projective curve
    of genus
    $g(B)$ with $s$ singular fibers.
    Assume that the smooth fibers $F$ are minimal, i.e., their
    canonical line bundles are semiample.
    Let $\omega_{X/B}$ be the relative canonical sheaf,
    and $\mathcal E\subseteq f_*\big(\omega_{X/B}^{\otimes k}\big)$
    be any non-zero subsheaf.
    Then the slope $\mu(\mathcal E):=\frac{\deg \mathcal E}{\rank
    \mathcal E}$ satisfies that
    $$\mu(\mathcal E)\leq \frac{mk\big(2g(B)-2+s\big)}{2}.$$
\end{theorem}

The main idea of proving \autoref{thm-high-dim} is to compute the
plurigenera by applying Riemann-Roch theorem for the direct image
sheaves $f_*\big(\omega_X^{\otimes k}\big)$ on the base curve.
Combining with the asymptotic behavior of the plurigenera, we
complete the proof.

\begin{proof}[Proof of \autoref{thm-high-dim}]
    Since the base is a rational curve $\bbp^1$,
    it follows that
    $$\omega_X=\omega_{X/\bbp^1}\otimes f^*\omega_{\bbp^1}
    =\omega_{X/\bbp^1}\otimes f^*\mathcal O_{\bbp^1}(-2).$$
    Hence for any $k\geq 1$, one has
    \begin{equation}\label{eqn-2-2}
        f_*\big(\omega_X^{\otimes k}\big)=
        f_*\big(\omega_{X/\bbp^1}^{\otimes k}\big)\otimes
        \mathcal O_{\bbp^1}(-2k).
    \end{equation}
    Let $\mathcal E\subseteq f_*\big(\omega_X^{\otimes k}\big)$
    be any subsheaf. Then by \eqref{eqn-2-2},
    $\mathcal E\otimes \mathcal O_{\bbp^1}(2k)$ is a
    subsheaf of $f_*\big(\omega_{X/\bbp^1}^{\otimes k}\big)$.
    Thus by \autoref{thm-ltz}, one obtains
    $$\mu(\mathcal E)+2k = \mu\big(\mathcal E\otimes
    \mathcal O_{\bbp^1}(2k)\big) \leq \frac{mk}{2}\cdot (s-2);$$
    equivalently, we have
    \begin{equation}\label{eqn-2-3}
        \mu(\mathcal E)\leq \frac{k}{2}\big(m(s-2)-4\big).
    \end{equation}

    As a locally free sheaf on $\bbp^1$, the direct image
    sheaf $f_*\big(\omega_X^{\otimes k}\big)$ is isomorphic to
    a direct sum of invertible sheaves,
    $$f_*\big(\omega_X^{\otimes k}\big)\cong \bigoplus_{i=1}^{r_k}
    \mathcal O_{\bbp^1}(d_i),
    \qquad  r_k=\rank f_*\big(\omega_X^{\otimes k}\big).$$
    By \eqref{eqn-2-3}, we have
    \begin{equation*}
    d_i\leq \frac{k}{2}\big(m(s-2)-4\big),\qquad
    i=1,2,\cdots, r_k.$$
    Hence
    $$\begin{aligned}
    &~\dim H^0(X,\omega_X^{\otimes k})=\dim H^0
    \big(\bbp^1,f^*(\omega_X^{\otimes k})\big)
    =\sum_{d_i\geq 0}(d_i+1)\\
    &~\leq \max\left\{0,~\Big[\frac{k}{2}
    \big(m(s-2)-4\big)\Big]+1\right\}\cdot
    \rank f_*\big(\omega_X^{\otimes k}\big).
    \end{aligned}
    \end{equation*}
    Here `$[\bullet]$' stands for the integral part.

    According to the definition of the Kodaira dimension of a variety,
    when $k$ is sufficiently large, one has
    $$\left\{\begin{aligned}
    &\rank f_*\big(\omega_X^{\otimes k}\big)
    =\dim H^0(F,\omega_F^{\otimes k})
    \sim k^{\kappa(F)};\\
    &\dim H^0(X,\omega_X^{\otimes k}) \sim k^{\kappa(X)}.
    \end{aligned}\right.$$
    Hence $\kappa(X)\leq \kappa(F)+1$.
    Moreover, if $\kappa(X)\geq 0$, then
    $$\frac{1}{2}\big(m(s-2)-4\big)\geq 0,\quad \text{~i.e.,~}\quad
    s\geq \frac{4}{m}+2;$$
    and if $\kappa(X)= \kappa(F)+1$, then
   $$\frac{1}{2}\big(m(s-2)-4\big)>0,\quad \text{~i.e.,~}\quad
   s>\frac{4}{m}+2.$$
   This completes the proof.
\end{proof}

\begin{remark}
    Recall that the volume of a projective variety $X$ is
    defined to be
    $${\rm Vol}(X)=\limsup_{k}~ \frac{(\dim X)! \cdot
    \dim H^0(X,\omega_X^{\otimes k})}{k^{\dim X}}.$$
    The above proof shows also that for a variety of
    general type semi-stably fibred over $\bbp^1$ with
    $s$ singular fibers,
    one has
    $${\rm Vol}(X)\leq \frac{(m+1)\big(m(s-2)-4\big)}{2}
    {\rm Vol}(F),$$
    where $F$ is a general fiber of $f$.
    In particular, when $X$ is of general type and
    $f:\,X \to \bbp^1$ is a semi-stable fibration of
    curves of genus $g\geq 2$ with $s$ singular fibers,
    one computes that
    $${\rm Vol}(X)=\omega_{X_0}^2,\qquad {\rm Vol}(F)=2g-2,$$
    where $X_0$ is the minimal model of $X$.
    Hence the above proof shows that in this case,
    $$\omega_{X_0}^2\leq 2(s-6)(g-1).$$
\end{remark}


\providecommand{\bysame}{\leavevmode\hbox to3em{\hrulefill}\thinspace}
\providecommand{\MR}{\relax\ifhmode\unskip\space\fi MR }
\providecommand{\MRhref}[2]{%
     \href{http://www.ams.org/mathscinet-getitem?mr=#1}{#2}
}
\providecommand{\href}[2]{#2}

\end{document}